\input amssym.def
\input amssym
\documentstyle[11pt]{article}
\begin{document}
\newtheorem{thm}{Theorem}[section]
\newtheorem{defin}[thm]{Definition}
\newtheorem{lemma}[thm]{Lemma}
\newtheorem{propo}[thm]{Proposition}
\newtheorem{cor}[thm]{Corollary}
\newtheorem{conj}[thm]{Conjecture}

\centerline{\huge \bf Examples of Categorification}
\bigskip

\centerline{\parbox{56mm}{Louis Crane and David N. Yetter \\{\em Department of
Mathematics \\ Kansas State University \\ Manhattan, KS 66506}}}
\bigskip
\bigskip

\section{Introduction}

The suggestion was made in work of Crane and Frenkel \cite{cf} \cite{c}
 that the
inverse relation to the Grothendieck rig (fusion rig) construction
should shed light on the relation between topological quantum field theories
(TQFT's) in various dimensions, and, as well, should provide constructions
for TQFT's in dimension 4. 

It is the purpose of this note to consider several simple cases of this 
relation. 

\begin{defin} A {\bf rig} $R$ is a set equipped with a monoid structure,
$(R, +, 0)$ and a semi-group structure $(R, \cdot)$ satisfying, moreover, 

\begin{center}
\begin{tabular}{rl}
commutativity & $a + b = b + a$ \\
right-distributivity & $(a + b)\cdot c = a\cdot c + b\cdot c$ \\
left-distributivity & $a\cdot (b + c) = a\cdot b + a\cdot c$. \\
\end{tabular}
\end{center}

\noindent A rig is {\bf unital} if the multiplicative semi-group is a monoid
with unit $1$. A rig is {\bf of finite rank} if the additive monoid $(R,+,0)$ 
is finitely generated.  
\end{defin}

Observe that in analogy to the coproduct of rings---tensor product of rings 
over ${\Bbb Z}$---rigs admit a coproduct, ``tensor product over ${\Bbb N}$''.

The one other general property of rigs (of finite rank) which we will need 
in the sequel is 

\begin{propo}
If $R$ is a rig of finite rank, then $R$ has a unique minimal set of additive
generators.
\end{propo}

\noindent{\bf proof:} Suppose $S$ and $T$ are minimal sets of generators.
Then there exist expressions

\[ s = \sum_{t \in T} n^s_t t \]

\noindent and

\[ t = \sum_{s \in S} m^t_s s \]

\noindent Observe then that $[n^s_t]$ and $[m^t_s]$ are mutually inverse
matrices of natural numbers.  The only such pairs are pairs of inverse
permutation matrices. $\Box$

\begin{defin}
If $\cal C$ is an (essentially small) tensor category, that is an
abelian category equipped with a (bi-)exact monoidal
productnot necessarily with unit object), 
then the set of isomorphism classes of 
objects in $\cal C$ equipped
with the operations induced by direct sum and tensor product is called
the {\bf Grothendieck rig} of $\cal C$, and denoted $Groth({\cal C})$. 
\end{defin}

	Notice that if $\cal C$ has a unit object,
then $Groth({\cal C})$ is unital. Similarly if $\cal C$ is Artinian 
semi-simple, then $Groth({\cal C})$ is of finite rank.  

Observe that the universal property of direct sum imposes some constraints
on the structure of a Grothendieck rig, the most salient of which is
the condition that all 1-generator sub-additive-monoids are free.

\begin{defin} An {\bf abstract fusion-rule algebra} is a rig, in which
all 1-generator sub-additive-monoids are free. 
An {\bf abstract fusion-rule bialgebra}
is an abstract fusion rule algebra $A$ equipped with a co-operation
$\Delta:A\rightarrow A\otimes_{\Bbb N} A$ 
which is a rig homomorphism (or, equivalently, which
satisfies the usual compatibility relations for bialgebras). It is counital
if it is equipped with a co-operation $\epsilon:A\rightarrow {\Bbb N}$ which
is a rig homomorphism. In the case where $A$ is unital, we require
$\Delta$ and $\epsilon$ to preserve the unit.
\end{defin}

This latter, the notion of an abstract fusion-rule bialgebra, is important
because of the construction proposed by Crane and Frenkel \cite{cf}
of state-sum invariants of 4-manifolds using as initial data 
categorical analogues of Hopf-algebras, and the result of Crane and Yetter
\cite{cy.alg.tqfts} showing that any 4D TQFT with factorization
at corners has as part of its structure a formal ``bialgebra category.''

In order to explain the structure of these analogues to Hopf-algebras, 
and to fix the context in which our reversal of the Grothendieck rig
construction will take place, first fix an algebraically closed field
$K$. We will briefly mention the difficulties which non-algebraically
closed fields present later.  Let $VECT_K$ or simply $VECT$ denote the
category of finite-dimensional vector-spaces over $K$ with its usual
monoidal structure.

\begin{defin} Let $VECT_K-mod$ denote the 2-category of all (small) 
Artinian semi-simple 
$K$-linear categories, that is $K$-linear categories  equivalent 
to finite powers of $VECT_K$  with exact functors and natural transformations
as 1- and 2-arrows. We refer to the objects of this 2-category as
$VECT$-modules.  Let $\boxtimes$ denote a chosen  bifunctor
from $VECT_K-mod^2$ to $VECT_K-mod$ which selects an object
${\cal C}\boxtimes{\cal D}$ equipped with a functor from
${\cal C} \times {\cal D}$ exact in each variable separately and universal
among such.
\end{defin}

The existence of $\boxtimes$ has been shown in Yetter \cite{catlinalg}, 
and can be readily verified
by mimicking the construction of the tensor product of vector-spaces, 
except that instead of identifying objects as one identified elements, one
must adjoin an isomorphism.  It follows readily from the universality
properties that $\boxtimes$ makes $VECT_K-mod$ into
a monoidal bicategory (cf. Gordon/Power/Street \cite{gps} and
Kapranov/Voevodsky \cite{kv}), whose underlying bicategory is a
2-category.

In this setting if a $VECT$-module $\cal C$ is equipped with a monoidal
structure (with unit object), exact in each variable, we may regard this as
 being given
by (exact) functors $\otimes: {\cal C}\boxtimes{\cal C}\rightarrow
{\cal C}$ (and $I: VECT\rightarrow {\cal C}$, sending $K$ to the monoidal
identity object), equipped with natural transformation(s) $\alpha$ (and
$\rho$ and $\lambda$) satisfying the usual pentagon (and triangle)
coherence condition(s). Similarly the structural transformations of
exact monoidal functors between $VECT$-modules $\cal C$ and $\cal D$
 may be understood as
natural transformations between functors from ${\cal C}\boxtimes
{\cal C}$ to $\cal D$ (and $VECT$ to $\cal D$). We will refer to monoidal
categories of this sort as {\em (Artinian semi-simple) tensor 
categories.}\footnote{More generally, we 
advocate the use of ``tensor category'' to
refer to a monoidal abelian category with $\otimes$ exact in each variable.}
We call exact monoidal functors {\em tensor functors}.

It is now sensible to consider in this setting duals to the notions of tensor
category and tensor functor:

\begin{defin}
An Artinian semi-simple {\bf (counital) cotensor category} 
over $K$ is a $VECT$-module
$\cal C$ equipped with functors $\Delta:{\cal C}\rightarrow {\cal C}
\boxtimes {\cal C}$ (and $\epsilon:{\cal C}\rightarrow VECT$), together
with natural transformations $\beta$ (and $r$ and $l$) 
satisfying the obvious pentagon
(and triangle) relation(s). 
A {\bf strong cotensor functor} $F:{\cal C}\rightarrow
{\cal D}$ is a functor equipped with natural isomorphisms $F_{\sim}:
(F\boxtimes F)(\Delta)\rightarrow \Delta(F)$ 
(and $F^0:\epsilon \rightarrow \epsilon(F)$) 
satisfying coherence conditions formally dual
to those for strong monoidal functors.
\end{defin}

We can now succinctly define the categorical analogue of a bialgebra
as given by Crane and Frenkel \cite{cf} including the conditions on
the unit and counit functors which were omitted in \cite{cf}:

\begin{defin}
An Artinian semi-simple {\bf bitensor category} over $K$ is an
Artinian semi-simple category equipped with both a tensor category 
structure and a cotensor category structure for which the structure
functors for the cotensor category strucure are (strong)
tensor functors, the structure functors for the tensor
structure are cotensor functors, and the structural transformations
for the tensor functor and cotensor functor structures coincide
whenever their sources and targets coincide. It is {\bf biunital} if the
tensor structure is unital and the cotensor structure is counital.
\end{defin}

In particular, we fix notation for these ``compatibility transformations''
as follows:

\[ \Phi = \tilde{\Delta} = \otimes_{\sim}^{-1} \]

\[ \eta = \epsilon^0 = I_0^{-1} \]

\[ \tau = \tilde{\epsilon} = \otimes_0^{-1} \]

\[ \delta = \Delta^0 = I_{\sim}^{-1} \]

The other two ``coherence cubes'' of \cite{cf} (besides the pentagon
and dual pentagon) are simply the coherence condition for a tensor
functor and its dual.

Thus the complete structure of a biunital bitensor 
category $\cal C$ is given by 
four functors  $\otimes, I, \Delta,$ and $\epsilon$ and
ten natural isomorphisms $\alpha, \rho, \lambda, \beta, r, l, \Phi, \eta,
\tau,$ and $\delta$ satisfying coherence conditions which can be read off
from the definitions of tensor and cotensor categories and tensor and
cotensor functors.

The appropriate notion of structure preserving functors between 
bitensor categories is given by

\begin{defin}
A {\bf (strong) bitensor functor} is a 5-tuple 
$(F,\tilde{F},F_{\sim},F^0,F_0)$, where $F:{\cal C}\rightarrow {\cal D}$
is a functor between bitensor categories, and
such that $(F,\tilde{F},F_0)$ (resp.
$(F,F_{\sim},F^0)$) is a tensor (resp. cotensor) functor and which moreover
satisfies the following, in which primes ($^\prime$) indicate structural
functors or natural transformations belonging to  $\cal D$:

\[ [ F\boxtimes F (\Phi_{A,B})][ \oplus \tilde{F}_{A_{(1)},A_{(2)}}\boxtimes
 \tilde{F}_{B_{(1)},B_{(2)}}] [ F_{\sim A}\boxtimes F_{\sim B}] = 
F_{\sim A\otimes B} \; \Delta^\prime(\tilde{F}_{A,B}) \; 
\Phi^\prime_{F(A),F(B)} \]

\[ F^0_{A\otimes B} \epsilon^\prime (\tilde{F}_{A,B}) \tau^\prime_{F(A),F(B)} 
= \tau_{A,B}
(F^0_A\otimes F^0_B) \]

\[ F_{\sim I} \Delta^\prime (F_0) \delta^\prime
 = (F\boxtimes F)(\delta) (F_0\boxtimes F_0) \]

\noindent and

\[ F^0_n F_0 \eta^\prime = \eta \]

If either category
is not unital or counital, the appropriate functor and attendant
natural isomorphisms are omitted. 
\end{defin}

Likewise, we can define cotensor and bitensor natural transformations:
a natural transformation is a cotensor transformation if it satisfies 
the dual condition to monoidal naturality, and is a bitensor transformation
if it is both a monoidal natural transformation and a cotensor transformation.

Now, observe that if $\cal C$ is a bitensor category, 
$Groth({\cal C})$ has the
structure of an abstract fusion rule bialgebra with the co-operations
induced by the cotensor and counit functors on the category. We will refer
to this as the {\bf Grothendieck birig} of the bialgebra category.

\begin{defin} A {\bf $K$-categorification} of an abstract fusion-rule
algebra (resp. bialgebra) $A$ is a $K$-linear tensor category
(resp. bialgebra category) whose Grothendieck rig (resp. birig) is $A$.
If $A$ is of finite rank, we call a categorification {\bf semi-simple} if it
is semi-simple as a tensor category, and, moreover, the objects whose images
under the Grothendieck rig construction are the additive generators of
$A$ are simple objects.
\end{defin}

Notice in general, 
in the non-finite rank case, one must specify a set of additive
generators to make sense of the notion of semi-simple categorification. We
will have no need of this more general notion in this paper.

In the remainder of the paper, we consider examples of categorifications
and their relevance to the construction of TQFT's. We restrict our 
attention to the case of $K$ algebraically closed because this
covers the most interesting case of $K = {\Bbb C}$ and removes the
possibility of having objects in semi-simple categories whose endomorphism
algebras are division-algebra extensions of the ground field.

\section{Categorifying ${\Bbb N}[G]$, Dijkgraaf-Witten Theory, and the
Turaev-Viro Construction}

Our first example sheds light on the relationship between two well-known
constructions of (2+1)-dimensional TQFT's:  Dijkgraaf-Witten theory
\cite{dw}, in particular its simplicial construction as in Wakui \cite{wak},
(cf. also Yetter \cite{yetgrps}), and the generalized Turaev-Viro
construction, cf. Barrett and Westbury \cite{bw}. As examples of
tensor categories, these are reasonably well-known (cf.  \cite{cp}, \cite{dpr})
and are included here merely as the simplest examples of categorification.

Fix a finite group $G$.  Observe that the group rig ${\Bbb N}[G]$ is 
an abstract fusion-rule bialgebra with the operations induced on the basis
$G$ by $g\cdot h = gh$ (the null-infix denoting the group-law) and
$\Delta(g) = g\otimes g$. The following shows that the 
categorifications of ${\Bbb N}[G]$ are essentially classified by the
3-cocyles on $G$ with coefficients in $K^\times$ when 
$K$ is algebraically closed.

	In this and all subsequent proofs, it first should be observed
that each equivalence class of $K$-linear semi-simple categories contains
skeletal categories (i.e. categories with only one object in each
isomorphism class).  We begin by restricting our attention to the structure of
these skeletal categories up to isomorphism, then consider the question
of monoidal equivalence.  Notice, we are taking an approach orthogonal
to that of the Mac Lane Coherence Theorem \cite{mac}:  we will not
be able to strictify the structure maps and retain skeletalness.

\begin{thm}
Let $G$ be a finite group, and $K$ an algebraically closed field, then the
isomorphism classes of skeletal semi-simple
$K$-categorifications of ${\Bbb N}[G]$ as an abstract fusion-rule 
algebra (resp. unital abstract fusion-rule algebra; 
abstract fusion-rule bialgebra; 
biunital abstract fusion-rule bialgebra) 
are in one-to-one correspondence with 
the set of $K^\times$-valued 3-cocycles on $G$ (resp.  
pairs $(\alpha, \rho)$, where $\alpha$ is a $K^\times$-valued 3-cocycle
on $G$ , and $\rho$ is an element of $K^\times$; $K^\times$-valued 2-cochains
on $G$; triples $(\phi, \rho, r)$, where $\phi$ is a $K^\times$-valued 
2-cochain on $G$, $\rho$ is an element of $K^\times$, and $r$ is a 
$K^\times$-valued 1-cochain on $G$).
\end{thm}

\noindent {\bf proof:} We begin by proving the first two statements, that
skeletal categorifications as an algebra (resp. unital algebra) are given by 
3-cocycles (resp. pairs of a 3-cocycle and a scalar).

	Now, in any skeletal categorification of ${\Bbb N}[G]$, we 
may identify the object whose image under the Grothendieck rig
construction is $g \in G$ with $g$. Since the tensor product (resp. 
identity object) must be carried to the multiplication (resp. unit)
in the fusion ring, the functors of the monoidal structure are given
by 

\[ g \otimes h = gh \;\;\; I = e \]

	To specify a categorification (as an algebra), it remains only 
to describe the rest of the monoidal structure:  in this case the
structure maps become families of maps $\alpha_{g,h,k}:ghk\rightarrow ghk$,
$\rho_g:g\rightarrow g$ and $\lambda_g:g\rightarrow g$. Observe, moreover,
that these ``maps'' are just elements of $K$, that
the semi-simplicity condition implies that any such families of
maps will satisfy the required naturality conditions, and that invertibility
consists in restricting the choices to elements of $K^\times$.

The pentagon condition for a 4-tuple $g,h,k,l$ can then be written:

\[ \alpha_{g,h,k} \alpha_{g,hk,l} \alpha_{h,k,l} = 
\alpha_{gh,k,l} \alpha_{g,h,kl} \]

\noindent which is precisely the condition that $\alpha_{-,-,-}$ be a
$K^\times$-valued 3-cocycle on $G$.

Upon including a unit, the triangle condition relating 
$\rho, \lambda,$ and $\alpha$ becomes

\[ \rho_g = \alpha_{g,e,h} \lambda_h \]

\noindent from which it follows that all components of $\rho$ and $\lambda$
are completely determined by $\alpha$ and the choice of a number 
$\rho$ such that $\rho_e = \rho \cdot Id_e$.  
The exercise of verifying that the choice of $\lambda_h$'s
determined by $\rho_e$, and the choice of $\rho_g$'s determined by
$\lambda_e$ satisfy the triangle condition for all pairs $g$ and $h$ is 
left to the reader.
(Hint: use the cocycle condition with two indices equal to $e$.)

	To specify a categorification as a bialgebra, notice first that
the dual pentagon condition on the coassociator reduces to

\[ \beta_g^3 = \beta_g^2 \]

\noindent for all $g \in G$, and thus since $\beta_g$ must be invertible,
$\beta_g = 1$.

	Similarly the compatibility condition between the connecting
transformation $\phi$ and the coassociator $\beta$ gives no restriction
on the components of $\phi$.

	On the other hand, the compatibility between $\phi$ and $\alpha$ 
is given by

\[ \phi_{g,h} \phi_{gh,k} \alpha_{g,h,k} = \alpha_{g,h,k} \alpha_{g,h,k}
\phi_{h,k} \phi_{g,hk} \]

\noindent which reduces to 

\[ \phi_{g,h} \phi_{g,hk}^{-1} \phi_{gh,k} \phi_{h,k}^{-1} = \alpha_{g,h,k} \]

	Thus, in this case the structure is completely determined by
the connecting transformation $\phi$, which is simply a 2-cochain on $G$.

	Finally in the case of a biunital bialgebra categorification,
the remaining structure maps are similarly completely determined by the
component maps at at simple objects, which are again given by scalars 
in $K^\times$,  
$r_g, l_g, \tau_{g,h}$; and since
 $e = I$
is simple, scalars $\delta$ and $\eta$.  The coherence conditions
then become:

\bigskip
\begin{center}
\begin{tabular}[h]{ll}
triangle for tensor structure & $\rho_g = 
\alpha_{g,e,k}\lambda_k$ \\
triangle for cotensor structure & $r_g = 
\beta_g l_g$ \\
$\Delta$ respects right unit & 
$\phi_{g,e}\delta\rho_g^2
 = \rho_g $ \\
$\Delta$ respects left unit & 
$\phi_{g,e}\delta\lambda_g^2
 = \lambda_g $ \\
$\otimes$ respects right counit & $
\phi_{g,k}\tau_{g,k} r_g
r_k = r_{gk} $\\
$\otimes$ respects left counit & $
\phi_{g,k}\tau_{g,k} l_g
l_k = l_{gk} $\\
$\epsilon$ preserves $\otimes$ & $ \alpha_{g,k,m}\tau_{g,km}\tau_{k,m} =
\tau_{g,k}\tau_{gk,m} $\\
$I$ preserves $\Delta$ & $ 
\beta_e \delta^2 
= \delta^2 $ \\
$\epsilon$ respects right unit & $\tau_{g,e}\eta = \rho_g$ \\
$\epsilon$ respects left unit & $\tau_{e,g}\eta = \lambda_g$ \\
$I$ respects right counit & $\delta \eta = r_e $\\
$I$ respects left counit & $\delta \eta = l_e $\\
\end{tabular}
\end{center}
\bigskip

We can then analyse these equations to determine a minimal set of
data and conditions for specifying the unital and counital structures
on a categorification of ${\Bbb N}[G]$. We assume that $\alpha$ and
$\phi$  have been chosen so as 
to specify a bitensor categorification without unit or counit.

First, observe that it follows from the two triangle conditions that
$\rho, \lambda, r$ and $l$ are completely determined by the values
of $\rho = \rho_e$  and $r_g$ $(g \in G)$ by
the formulas

\[ \rho_g = \alpha_{g,e,e}\rho \]

\[ \lambda_k = \alpha^{-1}_{e,e,k}\rho \]

\[ l_g = r_g \]

\noindent The only restriction on these needed to ensure that they
unambiguously determine cochains satisfying the first two equations
has already been imposed by the condition that $\alpha$  be a
coboundary in the suitable sense.

In a similar way,  $\otimes$ respecting
both the right and left counit conditions is equivalent to

\[ \tau(g,k) = d(r^{-1})(g,k) \phi^{-1}_{g,k}. \]

Similarly,  $\Delta$ respecting both the right and left
unit conditions is equivalent to

\[ \delta = \rho^{-1} \phi^{-1}_{e,e} \]

Given these last two equations, the conditions that $\epsilon$ preserve
$\otimes$ and $I$ preserve $\Delta$ follow from the 
 cocycle conditions and the fact that 
$d$ distributes over multiplication of cochains and (remember we write
the operation on cochains multiplicatively since the coefficients 
are in $K^\times$).

The four remaining conditions are all equivalent to 

\[ \eta = \rho r_e \phi_{e,e} \]

\noindent and we are done.$\Box$
\smallskip

	Turning to the more interesting question of monoidal equivalence
classes of categorifications, we have:

\begin{defin}
Two algebra (resp. bialgebra) categorifications are {\bf equivalent} if
there exists a monoidal equivalence (resp. bitensor 
equivalence) between them which induces the identity on Grothendieck
rigs.
\end{defin}

\begin{thm}
The equivalence classes of categorifications of ${\Bbb N}[G]$ as
an abstract fusion-rule algebra (whether unital or not) 
are in 1-1 correspondence with 
$H^3(G,K^\times)$.  The bialgebra categorification of ${\Bbb N}[G]$ is 
unique up to equivalence.
\end{thm}

\noindent {\bf proof:} Now, observe that the structure of a monoidal 
equivalence which induces the identity on the Grothendieck rig 
is given entirely by the structure transformations of the monoidal
equivalence.  In this case, these are given by a choice of units
$\psi_{a,b}$ for pair of group elements, and satisfying

\[ \alpha_{g,h,k}\psi_{h,k}\psi_{g,hk} = 
\psi_{g,h}\psi_{gh,k}\alpha^\prime_{g,h,k} \]

\noindent and

\[ \rho_g\psi_{g,e} = \rho^\prime_g \]

\noindent where $\alpha_{g,h,k}$ and $\alpha^\prime_{g,h,k}$ 
(resp. $\rho_g$ and $\rho^\prime_g$)\
are the components of
the associativity (resp. unit) transformations on the two categorifications.

Solving gives the condition that $\alpha$ and $a$ be cohomologous, while
the condition on the unit transformations is just a normalization condition
which can be trivially satisfied. 

The statement for bialgebra categorifications follows immediately from this
and the condition that the associativity constraint be the coboundary of the
connecting constraint:  $\psi$ must be the ratio of the two connecting
constraints. 

	The additional structure of unital-counital categorifications
is effaced by equivalence.

	A similar analysis to that
in the proof of the previous theorem shows that the structural transformations
for an equivalence of unital and counital categorifications are determined
by the 2-cochain $\psi$ which determines the non-unital non-counital
equivalence, together with a scalar $f_0$ and a 1-cochain $f^0$.

	The conditions besides the cobounding conditions 
 reduce to

\[ f_0 = 
\frac{\rho}{\rho^\prime\psi_{e,e}} \]

\[ f^0(g) = \frac{r_g}{r^\prime_g}
 \]

\noindent which can always be solved for any choice of $\rho$ and $\rho^\prime$
(resp. $r$ and $r^\prime$).
$\Box$

The algebra case of this result in fact shows that the Dijkgraaf-Witten 
invariants of 3-manifolds are examples of the generalized Turaev/Viro
construction of Barrett and Westbury \cite{bw}.  It is not hard to show
that the categorifications of a finite group rig are spherical
categories in the sense of \cite{bw}.   The construction as given by 
Wakui \cite{wak} is then immediately seen to be a special case of their 
general construction.

\section{Categorifying $D({\Bbb N}[G])$}

We now turn to the question of categorifying the simplest really
non-trivial Hopf algebra:  the Drinfel'd double of a finite (non-commutative)
group algebra (here taken with ${\Bbb N}$-coefficients to produce
a birig).

The same techniques may be used to classify the semi-simple
categorifications of arbitary bicrossproducts of a finite group algebra with
a dual finite group algebra (cf. Majid \cite{maj}).

Recall that the Drinfel'd double of a finite group algebra can
be constructed by taking as a basis pairs $(g,\hat{h})$, where
$g$ and $h$ are elements of the group, and $\hat{h}$ indicates the
element in the dual basis corresponding to $h$ as a basis element
in ${\Bbb N}[G]$, with structure maps given by

\begin{center}
\begin{tabular}{rl}
multiplication & $(g,\hat{h})\cdot (k,\hat{l}) = \delta_{k^{-1}hk,l}
(gk, \hat{l})$ \\
unit & $1 = \sum_h (e,\hat{h})$ \\
comultiplication & $\Delta (g,\hat{h}) = \sum_{kl = h} (g,\hat{k})\otimes
(g,\hat{l})$ \\
counit & $\epsilon (g,\hat{h}) = \delta_{h,e}$ \\
\end{tabular}
\end{center}

To categorify this as an abstract fusion-algebra (resp. -bialgebra) we
must first consider what are the source and target data for a component of
the associator (resp. components of the associator, coassociator and
compatiblity transformation).

The typical component of an associator is a map

\[ \alpha_{(g,\hat{h}),(k,\hat{l}),(m,\hat{n})}:((g,\hat{h})\otimes 
(k,\hat{l}))\otimes (m,\hat{n}) \rightarrow (g,\hat{h})\otimes ((k,\hat{l})
\otimes (m,\hat{n})) \]

\noindent Observe that the source and target objects are both $0$ unless
$l = k^{-1}hk$ and $n = m^{-1}lm$, in which case both are the object
$(gkm,\hat{n})$.

Thus the associator may be regarded as a family of non-zero scalars
$\alpha(g,k,m;\hat{n})$ giving the non-zero components as multiples
of the identity on $(gkm,\hat{n})$. (Note: a choice
of $g,k,m,$ and $\hat{n}$ contains enough information to recover
the source and target data for a non-zero component of the associator.)

Similarly, the coassociator has components given by maps (in the
tensor cube of the category) 

\[ \beta_{(g,\hat{h})}: \oplus_{pk = h} \oplus_{ij = p} (g,\hat{i})\boxtimes 
(g,\hat{j})\boxtimes (g,\hat{k}) \rightarrow \oplus_{iq = h} \oplus_{jk = q}
(g,\hat{i})\boxtimes (g,\hat{j})\boxtimes (g,\hat{k}) \]

Each such map is determined by its components on the various

\[ (g,\hat{i})\boxtimes (g,\hat{j})\boxtimes (g,\hat{k}) \]

\noindent and thus by a family of scalars $\beta(g; \hat{i}, \hat{j}, \hat{k})$
(As above, the given indices contain enough information to recover to which 
summand of which component of the coassociator this scalar belongs.)

Finally, the compatibility transformation or ``coherer'' 
has as components maps (in the
tensor square of the category)

\[ \phi_{(g,\hat{h}),(k,\hat{l})}: \oplus_{mn = l} (gk,\hat{m})\boxtimes 
(gk,\hat{n})
\rightarrow \oplus_{ab = h} (gk, \widehat{k^{-1}ak})\boxtimes
(gk, \widehat{k^{-1}bk})
\]
 
\noindent In this case, the sources and targets of components are 
given only in the case
where the source and target are non-zero, which happens precisely when
$k^{-1}hk = l$.

Thus the coherer is determined by a family of scalars 
$\phi(g,k;\hat{m},\hat{n})$.

Simply writing down the pentagon coherence condition on the associator
in terms of the scalars $\alpha(g,k,m;\hat{n})0$ gives

\[ \alpha(k,m,p;\hat{q}) \alpha(g,km,p;\hat{q}) 
\alpha(g,k,m;\widehat{pqp^{-1}}) = \alpha(gk,m,p;\hat{q}) 
\alpha(g,k,mp;\hat{q}) \]

Similarly, the dual pentagon in terms of the scalars becomes

\[ \beta(g;\hat{j},\hat{k},\hat{l}) \beta(g;\hat{i},\hat{jk},\hat{l})
\beta(g;\hat{i},\hat{j},\hat{k}) = \beta(g;\hat{ij},\hat{k},\hat{l})
\beta(g;\hat{i},\hat{j},\hat{kl}) \]

\noindent (A choice of 3-cocycle for each group element.)

The coherence condition provided by the compatibility transformation
as the structure transformation for $\Delta$ as a monoidal
functor becomes:

\[ \alpha(g,k,m;\hat{p})\alpha(g,k,m;\hat{q})\phi(k,m;\hat{p},\hat{q})
\phi(g,km;\hat{p},\hat{q}) = \phi(g,k;\widehat{mpm^{-1}},\widehat{mqm^{-1}})
\phi(gk,m;\hat{p},\hat{q})\alpha(g,k,m;\widehat{pq}) \]

(The two occurences of $\alpha$ on the left come from the associator
for ${\cal C}\boxtimes{\cal C}$.)

Similarly, the coherence condition provided by the compatibility transformation
as the structure transformation for $\otimes$ as a cotensor functor
becomes

\[ \phi(g,k;\hat{p};\hat{r})\phi(g,k;\widehat{pr},\hat{s})
\beta(gk;\hat{p},\hat{r},\hat{s}) = 
\beta(g;\widehat{kpk^{-1}},\widehat{krk^{-1}},\widehat{ksk^{-1}})
\beta(k;\hat{p},\hat{r},\hat{s})\phi(g,k;\hat{r},\hat{s})
\phi(g,k;\hat{p},\widehat{rs}) \]

These conditions and those involving the other structural transformations
become more intelligible if we introduce a general setting for such scalar
valued functions:

For any finite group, $G$, let $\hat{G}$ denote the set of characteristic
functions of 1-element subsets of $G$.  Let $C_{n,m}(G,K^\times)$ (or
$C_{n,m}$ when $G$ and $K$ are clear from context) denote the
abelian group of all functions from $G^n\times \hat{G}^m$ to $K^\times$.
The groups $C_{n,m}$ then form a double complex
(written multiplicatively)  with differentials 
$d_2:C_{n,m}\rightarrow C_{n,m+1}$, given by the Hochschild coboundary
in the ``hatted indices'', and $\tilde{d}_1:C_{n,m}\rightarrow C_{n+1,m}$
given by the same formula as Hochschild coboundary in the ``unhatted indices''
except that when the last index is dropped, all hatted indices are
left-conjugated by the dropped index.

In terms of these coboundary operations, the coherence conditions already
interpreted for categorifications of $D({\Bbb N}[G])$ become

\[ \tilde{d}_1(\alpha ) = 1 \]

\[ d_2(\beta) = 1 \]

\[ d_2(\alpha) = \tilde{d}_1(\phi) \]

\[ d_2(\phi) = \tilde{d}_1(\beta) \]

\noindent that is, the triple $(\alpha, \phi, \beta)$ forms a coboundary in
the total complex of the double complex 
$(C_{n,m},\tilde{d}_1,d_2$ $n,m \geq 1)$. We will index the cohomology of
the total complex of $(C_{n,m},...)$ by $n+m-1$, and denote the groups
by ${\Bbb H}_\bullet (G,K^\times)$

By identical methods, one can verify that the structure transformations
for the preservation of the tensor product and cotensor product by a 
bitensor functor with the identity as underlying functor are determined
by families of scalars $\tilde{f}(g,k;\hat{l})$ and $f_\sim(g;\hat{m},\hat{n})$
such that the ratios of corresponding structure maps satisfy

\[ \alpha \alpha^{\prime -1} = \tilde{d}_1 (\tilde{f}) \]

\[ \phi \phi^{\prime -1} = d_2 (\tilde{f}) \tilde{d}_1(f_\sim) \]

\noindent and

\[ \beta \beta^{\prime -1} = d_2(f_\sim). \]

Thus we have

\begin{thm} \label{cohom.thm}
The skeletal semi-simple $K$-categorifications 
of $D({\Bbb N}[G])$ as an abstract fusion-rule
algebra (resp. abstract fusion-rule bialgebra) 
are in one-to-one correspondence with the 3-coboundaries in
$(C_{n,1},\tilde{d}_1)$ (resp. the 3-coboundaries in the total
complex of the double complex $(C_{n,m},\tilde{d}_1,d_2)$). Moreover,
the equivalence classes of categorifications are in natural one-to-one
correspondence with the elements of the cohomology group $H_{3,1}$ 
(resp. ${\Bbb H}_3(G,K^\times)$).
\end{thm}

As was done explicitly for the associator, coherer, and coassociator above,
we can examine the components of each of the other structural transformations
at a simple object.  In this way we find that $\rho$ (resp. $\lambda$,
$r$, $l$) is determined by a $1,1$-cochain $\rho(g;\hat{h})$ (resp. 
$\lambda(g;\hat{h})  $, $r(g;\hat{h})  $, $l(g;\hat{h})  $). 
For example, the typical component of $\rho$ at a simple object is a map
from $(g;\hat{h})\otimes \oplus_k (e;\hat{k})$ to $(g;\hat{h})$, but
the source is just $(g;\hat{h})$, so the map is a scalar multiple of
$1_{(g;\hat{h})}$.

Likewise, $\delta$ is determined by $0,2$-cochain $\delta(\hat{k},\hat{l})$,
$\tau$ by a $2,0$-cochain, and $\eta$ by a single element of $K^\times$.
(Note: although the double complex actually used in defining categorifications
does not have a 0-row or 0-column, it is helpful here and in what 
follows to consider the larger complex which does, since much of what
is needed to handle the unital and counital structures is conveniently
phrased in terms of cochains and coboundaries in the larger complex.)

Again, by way of example, $\delta$ is a map from $\Delta(I)$ to $I\boxtimes I$
(the latter being the unit object in ${\cal C}\boxtimes {\cal C}$), that is
a map from $\oplus_{h \; kl=h} (e;\hat{k})\boxtimes (e;\hat{l})$ to
$\oplus_{k \; l} (e;\hat{k})\boxtimes (e;\hat{l})$, and is thus determined by
a choice of scalar for each summand (the simple summands of the two sides
being the same, and each occurring with multiplicity one), that is a
$0,2$-cochain.

Of course, these functions satisfy conditions equivalent to the coherence
conditions for the natural transformations they define.  It is an easy exercise
to write out each of the coherence conditions in turn, instantiate the
objects with simple objects (or simple summands of $I$, as appropriate), and
write out the corresponding equation on the cochains.

The table below summarizes the resulting equations:
\bigskip

\begin{center}
\begin{tabular}[h]{ll}
triangle for tensor structure & $\rho(g;\widehat{klk^{-1}}) = 
\alpha(g,e,k;\hat{l})\lambda(k;\hat{l})$ \\
triangle for cotensor structure & $r(g;\hat{k}) = 
\beta(g;\hat{k},\hat{e},\hat{l}) l(g;\hat{l})$ \\
$\Delta$ respects right unit & 
$\phi(g,e;\hat{k},\hat{l})\delta(\hat{k},\hat{l})\rho(g;\hat{k})
\rho(g;\hat{l}) = \rho(g;\widehat{kl}) $ \\
$\Delta$ respects left unit & 
$\phi(e,g;\hat{k},\hat{l})\delta(\hat{k},\hat{l})\lambda(g;\hat{k})
\lambda(g;\hat{l}) = \lambda(g;\widehat{kl}) $ \\
$\otimes$ respects right counit & $
\phi(g,k;\hat{l},\hat{e})\tau(g,k) r(g;\widehat{klk^{-1}})
r(g;\hat{l}) = r(gk;\hat{l}) $\\
$\otimes$ respects left counit & $
\phi(g,k;\hat{l},\hat{e})\tau(g,k) l(g;\widehat{klk^{-1}})
l(g;\hat{l}) = l(gk;\hat{l}) $\\
$\epsilon$ preserves $\otimes$ & $ \alpha(g,k,m;\hat{e})\tau(g,km)\tau(k,m) =
\tau(g,k)\tau(gk,m) $\\
$I$ preserves $\Delta$ & $ 
\beta(e;\hat{k},\hat{l},\hat{m})\delta(\hat{k},\widehat{lm})
\delta(\hat{l},\hat{m}) 
= \delta(\hat{k},\hat{l})\delta(\widehat{kl},\hat{m}) $ \\
$\epsilon$ respects right unit & $\tau(g,e)\eta = \rho(g;\hat{e})$ \\
$\epsilon$ respects left unit & $\tau(e,g)\eta = \lambda(g;\hat{e})$ \\
$I$ respects right counit & $\delta(\hat{h},\hat{e})\eta = r(e;\hat{h}) $\\
$I$ respects left counit & $\delta(\hat{e},\hat{h})\eta = l(e;\hat{h}) $\\
\end{tabular}
\end{center}
\bigskip

By way of example, the condition that $\Delta$ be a monoidal functor
includes the condition that $\Delta$ respect the right unit transformation.
Written out this become the equation

\[ \Phi_{A,I} (Id \otimes_2 \delta) \rho_2 = \Delta(\rho) \]

\noindent where the subscripts $_2$ indicate the corresponding structure
in ${\cal C}\boxtimes {\cal C}$.  Now, $\rho_2 = \rho \boxtimes \rho$,
and if  $A = (g;\hat{h})$ then all of the sources and targets of the
maps in the equation are $\oplus_{kl = h} (g;\hat{k})\boxtimes (g;\hat{l})$,
and each map is thus determined by a scalar for each triple 
$(g;\hat{k},\hat{l})$.  Writing out a component of the equation above 
gives

\[ \phi(g,e;\hat{k},\hat{l})\delta(\hat{k},\hat{l})\rho(g;\hat{k})
\rho(g;\hat{l}) = \rho(g;\widehat{kl}) \]

\noindent(Notice: on the right hand side, we use the fact that
$\Delta$ preserves identities
and scalar multiples.) 

We can then analyse these equations to determine a minimal set of
data and conditions for specifying the unital and counital structures
on a categorification of $D({\Bbb N}[G])$. We assume that $\alpha,
\phi$ and $\beta$ have been chosen as in the previous theorem
to specify a bitensor categorification without unit or counit.

First, observe that it follows from the two triangle conditions that
$\rho, \lambda, r$ and $l$ are completely determined by the values
of $\rho(e;\hat{m})$ $(m \in G)$ and $r(g;\hat{e})$ $(g \in G)$ by
the formulas

\[ \rho(g;\hat{m}) = \alpha(g,e,e;\hat{m})\rho(e;\hat{m}) \]

\[ \lambda(k;\hat{l}) = \alpha^{-1}(e,e,k;\hat{l})\rho(e;\widehat{klk^{-1}}) \]

\[ r(g;\hat{k}) = \beta(g;\hat{k},\hat{e},\hat{e})r(g;\hat{e}) \]

\[ l(g;\hat{k}) = \beta(g;\hat{e},\hat{e},\hat{k})r(g;\hat{e}) \]

\noindent The only restriction on these needed to ensure that they
unambiguously determine cochains satisfying the first two equations
has already been imposed by the condition that $\alpha$ and $\beta$ be
coboundaries in the suitable sense.

In a similar way,  $\otimes$ respecting
both the right and left counit conditions is equivalent to

\[ \tau(g,k) = \tilde{d}_1(r^{-1})(g,k;\hat{e}) 
\phi^{-1}(g,k;\hat{e},\hat{e}) \]

\noindent with no further conditions imposed on $r$ or $\phi$. 

On the other hand,  $\Delta$ respecting both the right and left
unit conditions is equivalent to

\[ \delta(\hat{k},\hat{l}) = d_2(\rho^{-1})(e;\hat{k},\hat{l}) 
\phi^{-1}(e,e;\hat{k},\hat{l}) \]

\noindent together with the condition that $\delta$ be invariant under
simultaneous conjugation of both indices by elements of $G$.

Given these last two equations, the conditions that $\epsilon$ preserve
$\otimes$ and $I$ preserve $\Delta$ follow from the condition that
$(\alpha,\phi,\beta)$ be a cocycle and the fact that $\tilde{d}_1$ and
$d_2$ distribute over multiplication of cochains and
$\tilde{d}_1^2 = 1$ and $d_2^2 =1$ (remember we write
the operation on cochains multiplicatively since the coefficients 
are in $K^\times$).

The four remaining conditions are all equivalent to 

\[ \eta = \rho(e;\hat{e})r(e;\hat{e})\phi(e,e;\hat{e},\hat{e}) \]

	A similar analysis shows that the structural transformations
for an equivalence of unital and counital categorifications are determined
by the 1,2-cochain and 2,1-cochain which determine the non-unital non-counital
equivalence, together with a 0,1-cochain $f_0$ and a 1,0-cochain $f^0$.

	The conditions besides the cobounding conditions of Theorem 
\ref{cohom.thm} reduce to

\[ f_0(\hat{k}) = 
\frac{\rho(e;\hat{k})}{\rho^\prime(e;\hat{k})\tilde{f}(e,e;\hat{k})} \]

\[ f^0(g) = \frac{r(g,\hat{e})}{r^\prime(g,\hat{e})f_\sim(g;\hat{e},\hat{e})}
 \]

\noindent Thus, it follows that within any equivalence class of 
categorifications, the {\em only} constraint upon the choice of 
the functions $\rho(e;\hat{h})$ and $r(g;\hat{e})$ is the condition 
that $\delta(\hat{k},\hat{l})$ be invariant under simultaneous 
conjugation of both indices.  We have thus almost shown

\begin{thm}
Every skeletal biunital semi-simple bitensor $K$-categorification 
of $D({\Bbb N}[G])$ is determined by a choice 
of a 3-cocycle $(\alpha,\phi,\beta)$
in the total complex 
of the double complex $C_{i,j}$, together with a choice of
functions $\rho(e;\hat{k}):\hat{G}\rightarrow K^\times$ and
$r(g;\hat{e}):G\rightarrow K^\times$ subject to
the condition that $\phi(e,e;\hat{k},\hat{l})d_2(\rho)(e,\hat{k},\hat{l})$
be invariant under simultaneous conjugation of the hatted indices. Conversely,
every such choice determines such a categorification up to isomorphism.  
The equivalence classes of
unital counital bitensor categorifications of $D({\Bbb N}[G])$ are
in natural one-to-one correspondence with ${\Bbb H}_3(G,K^\times)$.
\end{thm}

\noindent{\bf proof:}  For the two statements, it remains only to
observe that the condition on $\delta$ is equivalent to the given 
condition on $\phi$ and $\rho$.  

	The final statement requires a little more work.  By the preceding 
remark, it suffices to show that every cohomology class admits a
representative for which the $\rho(e;\hat{k})$ can be chosen so the 
invariance condition holds.  Let $(\alpha,\phi,\beta)$ be an
arbitrary 3-cocycle. Now, consider the 2-cochain $(1,f)$, where
$1$ is the constant 2,1-cochain, and $f(g;\hat{k},\hat{l}) = 
\phi^{-1}(e,e;\hat{k},\hat{l})$.  Multiplying $(\alpha,\phi,\beta)$ by
the (total) coboundary of $(1,f)$ give a cohomologous 3-cocycle
such that $\phi(e,e;\hat{k},\hat{l}) = 1$.  Thus any constant
$\rho(e;\hat{k})$ suffices. $\Box$

	Now observe that for any group $G$ there is at least one solution
to the required equations:  if we choose {\em all} of the families of
scalars to be identically $1$, we obtain a solution.  We will refer to
this and any equivalent bitensor categorifications 
as {\bf trivial categorifications} of $D({\Bbb N}[G])$.

	Of course, it behooves us to exhibit a non-trivial bitensor 
categorification, since we as yet have no examples. The simplest family
of such may be described as follows:  let all of the families of scalars
be identically $1$ except $\beta(g;\hat{i},\hat{j},\hat{k})$. 
Observe that all conditions not involving $\beta$ are trivially satisfied,
and that the conditions involving $\beta$ then reduce to those defining
other quantities in terms of $\beta$ and the other scalars, and 
the conditions

\[ \beta(g;\widehat{hih^{-1}},\widehat{hjh^{-1}},\widehat{hkh^{-1}}) 
\beta(h;\hat{i},\hat{j},\hat{k}) = \beta(gh;\hat{i},\hat{j},\hat{k}) \]

\noindent and

\[ \beta(g;\hat{i},\hat{j},\hat{k})\beta(g;\hat{i},\hat{jk},\hat{l})
\beta(g;\hat{j},\hat{k},\hat{l}) = \beta(g;\hat{ij},\hat{k},\hat{l})
\beta(g;\hat{i},\hat{j},\hat{kl}) \]

	Thus, we may regard $\beta(g;-,-,-)$ as a function from $G$
to 3-cocycles on $G$ (written with hatted indices), satisfying the first
equation.  In particular, any group homomorphism from $G$ to the (abelian)
group of 3-cocyles invariant under simultaneous conjugation gives such
a function.  

	For a specific example, let $G$ be any group with $C_2$, the 
cyclic group of order 2, as a quotient (e.g. $G = {\frak S}_n$). Call an 
element of $G$ even when its image in $C_2$ is the identity; odd otherwise.
In this case $\beta$ given by

\begin{eqnarray*}
\beta(g;\hat{i},\hat{j},\hat{k}) & = & -1 \;\; 
\parbox{2in}{\rm if $g,i,j,k$ are all odd} \\
	& = & 1 \;\;\; \parbox{2in}{\rm otherwise} \\
\end{eqnarray*}

\noindent has all of the desired properties. In particular, $(1,1,\beta)$
represents a non-trivial element of ${\Bbb H}_3(G,K^\times)$ (provided 
$char(K) \neq 2$).

\section{Conclusions}

	The cohomological setting which provided a natural setting for these
constructions and classification theorems suggests that the process of
categorification should, at least in the semi-simple case, be viewed as
a deformation process for tensor or bitensor categories.

	The authors, in work in preparation \cite{cy.prep} have constructed
a general framework for the infinitesimal deformation of general 
(semi-simple) bitensor categories in terms of a similar 
double complex, and have isolated the cohomological obstructions to the
existence of formal power-series deformations as classes in the total
cohomology of the double complex.

	It still remains to use the examples presented herein to provide
explicit examples of $4$-manifold invariants of Crane/Frenkel type \cite{cf}
 and to 
construct the monoidal bicategory of representations of the bitensor
categories constructed herein, thereby giving initial data for a
fully bicategorical version of the Crane-Yetter construction (cf. \cite{cy},
\cite{cky}).

\end{document}